\documentclass[11pt,reqno,a4paper]{amsart}
\usepackage{color}
\usepackage{amssymb,amsmath}
\usepackage{bbm}
\usepackage[latin1]{inputenc}
\usepackage[active]{srcltx}
\usepackage{graphicx}
\usepackage{exscale,relsize}
\usepackage{textgreek}
\usepackage{epsfig,graphics}
\usepackage{psfrag}
\usepackage{caption}
\usepackage{subcaption}
\usepackage[toc,page]{appendix}
\usepackage{bigints}

\def\R{\mathbb{R}}
\def\m1{{I\!\!M}}

\def\ee{\`e}

\def\aa{\`a}

\newcommand{\grad}{\nabla}

\renewcommand{\to}{\rightarrow}
\newcommand{\pa}{\partial}

\newcommand{\ino}{\int_{\Omega}}

\newcommand{\rife}[1]{(\ref{#1})}
\newcommand{\ov}[1]{\overline{#1}}

\newcommand{\sscp}{\scriptscriptstyle}
\newcommand{\dsp}{\displaystyle}

\renewcommand{\dfrac}{\displaystyle\frac}
\newcommand{\finedim}{\hspace{\fill}$\square$}
\newcommand{\intbar}{\mathop{\int\makebox(-15.5,0){\rule[6pt]{.7em}{0.3pt}}\kern-6pt}\nolimits}

\newcommand{\ii}{\infty}

\newcommand{\al}{\alpha}

\newcommand{\ga}{\gamma}
\newcommand{\om}{\Omega}
\newcommand{\lm}{\lambda}

\newcommand{\rl}{\mbox{\Large \textrho}_{\!\sscp \lm}}

\newcommand{\thl}{\theta_{\sscp \lm}}

\renewcommand{\rho}{\mbox{\Large \textrho}}
\newcommand{\rh}{\mbox{\Large \textrho}}

\newcommand{\pl}{\psi_{\sscp \lm}}
\newcommand{\xil}{\ul}
\newcommand{\ul}{u_{\sscp \lm}}

\newcommand{\val}{v_{\sscp I}}

\newcommand{\ssl}{\sscp \lm}

\newcommand{\all}{\al_{\ssl}}
\newcommand{\gal}{\ga_{\sscp I}}
\newcommand{\el}{E_{\ssl}}

\newcommand{\fbi}{{\bf (F)$_{I}$}}
\newcommand{\prl}{{\textbf{(}\mathbf P\textbf{)}_{\mathbf \lm}}}

\newtheorem{theorem}{Theorem}[section]
\newtheorem{proposition}[theorem]{Proposition}
\newtheorem{lemma}[theorem]{Lemma}
\newtheorem{corollary}[theorem]{Corollary}
\newtheorem{remark}[theorem]{Remark}
\newtheorem{definition}[theorem]{Definition}
\newcommand{\brm}{\begin{remark}\rm}
\newcommand{\erm}{\end{remark}}
\newcommand{\bdf}{\begin{definition}\rm}
\newcommand{\edf}{\end{definition}}
\newcommand{\bte}{\begin{theorem}}
\newcommand{\ete}{\end{theorem}}
\newcommand{\bpr}{\begin{proposition}}
\newcommand{\epr}{\end{proposition}}
\newcommand{\ble}{\begin{lemma}}
\newcommand{\ele}{\end{lemma}}
\newcommand{\bco}{\begin{corollary}}
\newcommand{\eco}{\end{corollary}}
\newcommand{\beq}{\begin{equation}}
\newcommand{\eeq}{\end{equation}}
\newcommand{\bdm}{\begin{displaymath}}
\newcommand{\edm}{\end{displaymath}}

\newcommand{\graf}[1]{\left\{\begin{array}{ll}#1\end{array}\right.}

\def\sideremark#1{\ifvmode\leavevmode\fi\vadjust{\vbox to0pt{\vss
 \hbox to 0pt{\hskip\hsize\hskip1em \vbox{\hsize2.1cm\tiny\raggedright\pretolerance10000 \noindent #1\hfill}\hss}\vbox to15pt{\vfil}\vss}}}

\begin{document}
\numberwithin{equation}{section}
\parindent=0pt
\hfuzz=2pt
\frenchspacing

\title[Universal estimates for free boundary problems]{
New universal estimates for free boundary problems\\ arising in plasma physics}


\author[D. Bartolucci]{Daniele Bartolucci$^{(\dag)}$}
\address{Daniele Bartolucci, Department of Mathematics, University of Rome \emph{"Tor Vergata"}, Via della ricerca scientifica n.1, 00133 Roma.}
\email{bartoluc@mat.uniroma2.it}

\author[A. Jevnikar]{Aleks Jevnikar}
\address{Aleks Jevnikar, Department of Mathematics, Computer Science and Physics, University of Udine, Via delle Scienze 206, 33100 Udine, Italy.}
\email{aleks.jevnikar@uniud.it}

\thanks{$^{(\dag)}$Research partially supported by:
Beyond Borders project 2019 (sponsored by Univ. of Rome "Tor Vergata") "{\em Variational Approaches to PDE's}",
MIUR Excellence Department Project awarded to the Department of Mathematics, Univ. of Rome Tor Vergata, CUP E83C18000100006.}

\begin{abstract}
For $\om\subset \R^2$ a smooth and bounded domain, we derive a sharp universal
energy estimate for non-negative solutions of free boundary problems on $\om$
arising in plasma physics.
As a consequence, we are able to deduce new universal estimates for this class of problems.
We first come up with a sharp positivity threshold which guarantees that there is no free
boundary inside $\om$ or either, equivalently,
with a sharp necessary condition for the existence of a free boundary in the interior of $\om$.
Then we derive an explicit bound for the $L^{\infty}$-norm of non-negative solutions and also
obtain explicit estimates for the thresholds relative to other neat density boundary values.
At least to our knowledge, these are the first explicit estimates of this sort in the superlinear case.
\end{abstract}
\maketitle



\

\

\setcounter{section}{0}
\setcounter{equation}{0}
\section{\bf Introduction}

\

Letting $\om\subset \R^2$ be an open and bounded domain of class $C^{3}$,
we consider the free boundary problem
$$
\graf{-\Delta v = (v)_{+}^p\quad \mbox{in}\;\;\om\\ \\
-\bigintss\limits_{\pa\om} \dfrac{\pa v}{\pa\nu}=I \\ \\
v=\ga \quad \mbox{on}\;\;\pa\om
}\qquad \qquad \mbox{\bf (F)}_{I}
$$
for the unknowns $\ga \in \R$ and
$v\in C^{2,r}(\ov{\om}\,)$, $r\in (0,1)$. Here $(v)_+$ is the positive part of $v$,
$\nu$ is the exterior unit normal, $I> 0$ and $p\in (1,+\infty)$ are fixed.
Up to a suitable rescaling, we can assume without loss of generality that $|\om|=1$ and $(v)^p_+$ to be
multiplied by any positive constant.\\

The problem \fbi\, arises in Tokamak's plasma physics and we refer to \cite{Kad,Mer,Te} for a physical description of the problem.
A systematic analysis of \fbi\, has been initiated in \cite{BeBr,Te,Te2}. In particular, the authors in \cite{BeBr} considered the problem with more general operators and nonlinearities and showed that for any $I>0$ there exists at least one solution of \fbi. For old and new results about \fbi\, for
$p>1$, see for example \cite{AmbM,BMar,BMar1,BSp,FW,Mar,RenWei,We}, while for the model case $p=1$ (which requires a slightly different formulation, see the discussion after $\prl$)
\cite{dam,Gal,Pudam}. For further references and for the last developments about the uniqueness of solutions and about the qualitative behavior of the branch of solutions via bifurcation analysis, see \cite{BJ2}.\\

We will be here mainly concerned with positive solutions of \fbi,
which are related to the following dual formulation introduced in \cite{BeBr,Te2},
$$
\graf{-\Delta \psi =(\al+{\lm}\psi)^p\quad \mbox{in}\;\;\om\\ \\
\bigintss\limits_{\om} {\dsp \left(\al+{\lm}\psi\right)^p}=1\\ \\
\psi>0 \quad \mbox{in}\;\;\om, \quad \psi=0 \quad \mbox{on}\;\;\pa\om \\ \\
\al\geq0
}\qquad \prl
$$
for the unknowns $\al\in\R$ and $\psi \in C^{2,r}_{0,+}(\ov{\om}\,)$. Here, $\lm\geq 0$ and ${p\in [1,+\infty)}$ are fixed and for $r\in (0,1)$ we set
$$
C^{2,r}_0(\ov{\om}\,)=\{\psi \in C^{2,r}(\ov{\om}\,)\,:\, \psi=0\mbox{ on }\pa \om\},\;
C^{2,r}_{0,+}(\ov{\om}\,)=\{\psi \in C^{2,r}_0(\ov{\om}\,)\,:\, \psi> 0\mbox{ in } \om\}.
$$
Indeed, the relation between the dual problems \fbi\, and {\rm $\prl$} is as follows. Take $q$ such that
$$
\frac1p+\frac1q=1.
$$
For any fixed $\lm>0$ and $p>1$, $(\all,\pl)$ is a
solution of {\rm $\prl$} if and only if,
for $I=I_{\ssl}=\lm^{q}$, $(\gal,\val)=(\lm^{\frac{1}{p-1}}\all,\lm^{\frac{1}{p-1}}(\all+\lm\pl))$ is a
non-negative solution, i.e. with $\gal\geq0$,  of \fbi. Therefore
in particular, if $(\gal,\val)$ solves \fbi\, with $\gal \geq 0$, then
$(\all,\pl)=(I^{-\frac{1}{p}}\gal,I^{-1}(\val -\gal))$ solves
{\rm $\prl$} and the identity $I^{-\frac{1}{p}}\val=\all+\lm\pl$ holds. Finally, observe that for $p=1$ {\rm $\prl$} is already equivalent to a more general problem than \fbi\, and solutions of {\rm $\prl$} correspond to non-negative solutions of \fbi\, where the first equation is replaced by $-\Delta v = \lm(v)_{+}$.\\

We point out that since $|\om|=1$ and $\lm\geq 0$, then any solution $(\all,\pl)$ of $\prl$ satisfies
$$
\all\leq 1,
$$
and the equality holds if and only if $\lm=0$, for which $\prl$ admits a unique solution
which we denote by $\psi_{\sscp 0}$. The energy associated to a solution $(\all,\pl)$ of $\prl$ is defined as
$$
	\el=\frac12\int_{\om}|\grad \pl|^2.
$$

Here and in the rest of this paper $\mathbb{D}$ will denote the two-dimensional ball of unit area.
We will state the results in terms of $\prl$, keeping in mind the above discussed equivalence with \fbi.\\

Our first result is the following sharp universal energy estimate for any solution of $\prl$, depending only on the exponent $p$.
\bte\label{energy} Let $p\in [1,+\ii)$ and $(\all, \pl)$ be a solution of {\rm $\prl$}. Then it holds,
\beq\label{level0}
 2\lm\left(\frac{p+1}{16\pi}-\el\right)\geq \all(1-\all^p),
\eeq
where the equality holds if and only if, up to a translation, $\om=\mathbb{D}$. In particular,
\beq\label{level01}
\el\leq \frac{p+1}{16\pi},
\eeq
and the equality holds if and only if, up to a translation, $\om=\mathbb{D}$ and $\all=0$.
\ete

By making use of the latter result we will derive new universal estimates for this class of problems.
First of all the sharp character of \rife{level01} yields other sharp estimates for the positivity
threshold of solutions of $\prl$ and in particular of variational solutions of
$\prl$ and \fbi, as introduced in \cite{BeBr,Te2}, see also \cite{BMar,BMar1,BMar2,BSp}. For any plasma density
$$
\rh\in\mathcal{P}_{\sscp \om}:=\left\{\rho\in L^{1+\frac1p}(\om)\,|\,\rho\geq 0\;\mbox{a.e. in}\;\om \right\},
$$
and any $\lm\geq  0$, we define the free energy,
\beq\label{jeil}
J_{\ssl}(\rh)=
{\scriptstyle \frac{p}{p+1}}\ino (\rh)^{1+\frac{1}{p}}-\frac\lm 2 \ino \rho G[\rho],
\eeq
where $G[\rho](x)=\ino G_{\om}(x,y)\rho(y)\,dy$ and $G_{\om}$ is the Green function of $-\Delta$ with Dirichlet boundary conditions on $\om$. We then consider the minimization problem
$$
\mathcal{J}(\lm)=\inf\left\{ J_{\ssl}(\rh)\,:\,\rh\in \mathcal{P}_{\sscp \om}, \ino \rh=1\right\}.
$$
We know from \cite{BeBr,Te2} that for each $\lm>0$ there exists at least one $\rl$ which minimize $J_{\ssl}$.
In particular those minimizers $\rl\!\!$ whose Lagrange multiplier $\all$ is non negative yield a solution $(\all,\pl)$ of $\prl$ where $\pl=G[\rl]$. Any such solution is called a variational solution of $\prl$. Again, for $p>1$, there is an equivalent dual variational principle for \fbi\, which we will
not discuss here, see \cite{BeBr} and Appendix A in \cite{BJ2} for further details.\\

Concerning the positivity threshold for variational solutions we know by \cite{BMar1,Te2} (see also Corollary~A.1 in \cite{BJ2}) the following,

\medskip

{\bf Theorem A} (\cite{BMar1,Te2}). {\it
Let $p\in[1,+\infty)$ and $(\all,\pl)$ be a variational solution of {\rm $\prl$}.
Then there exists $\lm^{**}(\om,p)\in (0,+\ii)$ such that $\all>0$
if and only if $\lm\in (0,\lm^{**}(\om,p))$ and $\all=0$ if and only if $\lm=\lm^{**}(\om,p)$.}\\

Also, let $S_p(\om)$ be the best constant in the Sobolev embedding $\|w\|_p\leq S_p(\om)\|\nabla w\|_2$,
$w\in H^1_0(\om)$ and for $p\in[1,+\ii)$ let us define $\Lambda(\om,p)=S^{-2}_p(\om)$ and
$$
\lm^{2p}_{*}(\om,p)
=\left(\frac{8\pi}{p+1}\right)^{p-1}\Lambda^{p+1}(\om,p+1).
$$
It is well known (\cite{CRa1}) that $\lm_*(\om,p)\geq \lm_*(\mathbb{D},p)$ where
the equality holds if and only if,
up to a translation, $\om=\mathbb{D}$. Here we have,
\bte\label{thmvar} Let $p\in [1,+\ii)$
and $(\lm,\pl)$ be a solution of {\rm $\prl$} with $\lm\leq \lm_{*}(\om,p)$. Then
$\all>0$ unless either $p=1$ and $\lm=\lm_{*}(\om,1)$ or $p>1$, $\lm=\lm_*(\mathbb{D},p)$ and,
up to a translation, $\om=\mathbb{D}$, in which cases $\all=0$.\\ In particular, for variational solutions of {\rm $\prl$}, we have
$\lm^{**}(\om,p)\geq \lm_{*}(\om,p)$ and the equality holds if and only if either $p=1$ or $p>1$ and
up to a translation $\om=\mathbb{D}$.
\ete
It is interesting to comment about the sharp character of these estimates.
For $p=1$ they are sharp as they reduce to the well-known (\cite{BeBr,Te2}) sharp positivity threshold
$\lm^{**}(\om,1)=\lm_*(\om,1)=\lm^{(1)}(\om)\equiv \Lambda(\om,2)$, where $\lm^{(1)}$ is the first eigenvalue
of $-\Delta$ with Dirichlet boundary conditions. For $p>1$ they are sharp in the sense that
the equality $\all=0$ is attained
if and only if $\lm=\lm_*(\mathbb{D},p)$ and up to a translation $\om=\mathbb{D}$.
Moreover, since solutions of $\prl$/\fbi\, are unique on $\mathbb{D}$ (\cite{BSp}),
then interestingly enough, we obtain in this case the explicit value of the
sharp positivity threshold $\lm^{**}(\mathbb{D},p)=\lm_{*}(\mathbb{D},p)$, see also \cite{BJ2} for new results
concerning this case. Actually it seems that explicit estimates
about $\lm^{**}(\om,p)$ and $\lm_{*}(\om,p)$ were not known so far for $p>1$, while
besides the model case $p=1$, results of this sort are well known in the
sublinear case (\cite{AmbM}), see also \cite{BMar2} and \cite{CF}.
In particular, via the equivalence with \fbi, in the superlinear case $p>1$ we obtain a sharp
condition which guarantees that for variational solutions of \fbi\, there is no
free boundary inside $\om$ or either, equivalently, a necessary condition for the existence of a
free boundary in the interior of $\om$. Indeed, based on Theorems A and \ref{thmvar},
it is straightforward to deduce the following,
\bco
Let $p>1$ and $(\ga_{\sscp I}, v_{\sscp I})$ be a variational solution of
{\rm \fbi} with $\ga_{\sscp I}\leq 0$.
Then $I\geq (\lm_{*}(\om,p))^{\frac1q}$ and the equality holds
if and only if $\ga_{\sscp I}=0$ and up to a translation $\om=\mathbb{D}$.
\eco

\smallskip

Other useful information can be derived from Theorem \ref{energy}.
It is known by \cite{BeBr} that the solutions of \fbi\, and of $\prl$\,
are uniformly a priori bounded. However
such bound is obtained by standard elliptic estimates and bootstrap arguments
and thus explicit estimates were missing so far.
Our goal concerning this point
is to derive universal (independent on $\om$ and depending only on the exponent $p$)
explicit estimates for the $L^{\infty}$-norm of solutions of $\prl$.
On the other hand, one may ask which are the thresholds
relative to other neat values of $\all$. Let us introduce
\begin{equation} \label{l}
\ell(\om)=\frac{1}{2\pi}|\pa\om|^2-1,
\end{equation}
which, by the isoperimetric inequality, satisfies $\ell(\om)\geq \ell(\mathbb{D})=1$. Then we have,
\bte\label{corplasma1}
Let $p\in[1,+\infty)$ and let $(\all,\pl)$ be a solution of {\rm $\prl$}. Then the following holds.
\begin{itemize}
\item[1.] \emph{($L^{\infty}$ bound):}
\beq\label{est.intro}
 \|\pl\|_{\ii}<\frac{ p+1}{4\pi}\left(1+\frac{\lm p}{8\pi}\right).
\eeq
\item[2.] \emph{($\frac12$-$\frac1q$ thresholds):}
\begin{equation} \label{pos2}
\begin{array}{ll}
\all>\dfrac{1}{2} &\mbox{\emph{for} } \lm \leq \dfrac{4\pi}{e p}\,, \\
\all>\max\left\{\dfrac{1}{2},\dfrac{1}{q}\right\}  &\mbox{\emph{for} } \lm \leq \dfrac{4\pi}{e p\,\ell(\om)}\,.
\end{array}
\end{equation}
\end{itemize}
\ete

\bigskip

We point out that, concerning the $L^{\infty}$-bound \eqref{est.intro},
a slightly better estimate holds true, see Proposition~\ref{aprenrg}.
At least to our knowledge, these are the first explicit estimates
of this sort in the superlinear case $p>1$.
It is
interesting to check how far we get, with the argument pursued in Theorem \ref{corplasma1},
from the optimal result of Theorem \ref{thmvar}. Indeed, with this argument we can prove,
see Proposition \ref{lemc12},
that if $\lm \leq \frac{16\pi}{e (p+1)}$ then $\all>0$.

\bigskip

The main idea of the proof is to exploit the role of the energy $\el$ associated to a
solution $(\all,\pl)$ of $\prl$, which turns out to be related to the density
interaction energy  $\mathcal{E}(\rho)=\frac12\ino \rho G[\rho]$ for a
plasma density {\rm $\rh\in L^{1}(\om)$}.
Indeed, it is easily seen that $\el=\mathcal{E}(\rl)$ whenever $\pl=G[\rl]$.
First, we derive the sharp energy estimate, which is based on a differential inequality
involving level sets of solutions of $\prl$. This yields, by using also the isoperimetric
property of the Sobolev constant \cite{CRa1}, Theorem \ref{thmvar}.
Next, we manage to control the $L^{\infty}$-norm of a solution by means
of its energy and then the uniform bound follows once more by the energy estimate.
Then we exploit the $L^{\infty}$-bound to deduce the $\frac12$-$\frac1q$ thresholds.

\

\

\section{\bf Proof of the main results}\label{sec2}

\

We collect in this section the proof of the main results. We divide the argument in several steps. Letting $(\all, \pl)$ be a solution of {\rm $\prl$}, it will be convenient to use the auxiliary function $\xil=\lm \pl$ which satisfies
\beq\label{ul.1}
\graf{-\Delta \xil =\lm \left(\all+\xil\right)^p\quad \mbox{in}\;\;\om\\ \\
\bigintss\limits_{\om}  \left(\all+\xil\right)^p=1\\ \\
\xil\geq0 \;\; \mbox{in}\;\;\om, \quad \xil=0 \;\; \mbox{on}\;\;\pa\om \\ \\
\all\geq0.
}
\eeq
We point out that we already know by \cite{BeBr} that the $L^{\infty}$-norm of $\xil$ is uniformly bounded (see also Proposition \ref{aprenrg} below) and so we assume without loss of generality that
$$
\thl:=\xil(0)=\|\xil\|_{L^\ii(\om)}.
$$
Since $\ino \left(\all+{\xil}\right)^p=1$ and $|\om|=1$, then necessarily for any solution and for $\lm>0$,
\beq\label{nm1}
\left(\all+{\thl}\right)^p>1.
\eeq

We start with the sharp universal energy estimate.

\medskip

\proof[Proof of Theorem \ref{energy}]
Let $\thl=\|\xil\|_{L^{\infty}(\om)}$ and set
$$
\om(t)=\{x\in\om\,:\,\xil>t\}, \quad \Gamma(t)=\{x\in\om\,:\,\xil=t\}, \quad t\in [0, \thl],
$$
and
$$
m(t)=\lm \int\limits_{\om(t)}\left(\all+{\xil}\right)^{p},\qquad \mu(t)=|\om(t)|,
\qquad {e}(t)=\int\limits_{\om(t)}|\nabla \xil|^2,
$$
where $|\om(t)|$ is the area of $\om(t)$. If $\lm=0$ then \rife{level0} is trivially satisfied and \rife{level01} follows by well known rearrangement estimates (\cite{Tal2}),
\beq\label{stima-energy}
E_{0}(\om)=\frac12\ino \ino G_{\om}(x,y)\,dxdy\leq
\frac12\int_{\mathbb{D}} \int_{\mathbb{D}} G_{\mathbb{D}}(x,y)\,dxdy=
\frac{1}{16\pi}.
\eeq

Hence, we consider now $\lm>0$. Since $|\Delta\xil|$ is bounded below away from zero and since the boundary is smooth,
then it is not difficult to see that actually $m(t)$ and $\mu(t)$ are continuous in $[0, \thl]$
and piecewise smooth in $[0, \thl]$, that is, of class $C^{1}$ with the exception of a finite number of points in $[0,\thl]$.
In particular the level sets have vanishing area $|\Gamma(t)|=0$ for any $t$ and we will use the fact that,
$$
m(0)=\lm, \quad \mu(0)=1, \quad e(0)=\int\limits_{\om}|\nabla \xil|^2\equiv 2\lm^2\el,
$$
and
$$
m(\thl)=0, \quad \mu(\thl)=0, \quad e(\thl)=0.
$$
By the co-area formula and the Sard Lemma we have,
\beq\label{diffm}
-m^{'}(t)=\lm \int\limits_{\Gamma(t)}\frac{\left(\all+{\xil}\right)^{p}}{|\nabla \xil |}=
\lm \left(\all+{t}\right)^{p}\int\limits_{\Gamma(t)}\frac{1}{|\nabla \xil |}=
\lm \left(\all+{t}\right)^{p}(-\mu^{'}(t)),
\eeq
and
\beq\label{diffem}
m(t)=-\int\limits_{\om(t)}\Delta \xil=\int\limits_{\Gamma(t)}|\nabla \xil|=-e^{'}(t),
\eeq
for a.a. $t\in [0, \thl]$.
By the Schwarz inequality and the isoperimetric inequality we find that,
$$
-m^{'}(t)m(t)=
\lm \int\limits_{\Gamma(t)}\frac{\left(\all+{\xil}\right)^{p}}{|\nabla \xil |}\int\limits_{\Gamma(t)}|\nabla \xil|=
\lm \left(\all+{t}\right)^{p}\int\limits_{\Gamma(t)}\frac{1}{|\nabla \xil |}\int\limits_{\Gamma(t)}|\nabla \xil|\geq
$$
$$
\lm \left(\all+{t}\right)^{p}\left(|{\Gamma(t)}|_{1}\right)^2\geq\lm\left(\all+{t}\right)^{p}
4\pi \mu(t),
\mbox{ for a.a. } t\in [0, \thl],
$$
where $|{\Gamma(t)}|_{1}$ denotes the length of $\Gamma(t)$. Therefore, we conclude that,
\beq\label{diff1}
\frac{(m^2(t))^{'}}{8\pi}+\lm \left(\all+{t}\right)^{p}\mu(t)\leq 0,\mbox{ for a.a. } t\in [0, \thl].
\eeq
By using the following identity,
$$
\left(\all+{t}\right)^{p}\mu(t)=\frac{1}{p+1}\left(\left( \all+{t}\right)^{p+1}\mu(t)\right)^{'}
-\frac{1}{p+1}\left(\all+{t}\right)^{p+1}\mu^{'}(t), \mbox{ for a.a. } t\in [0, \thl],
$$
together with \rife{diff1} and \rife{diffm} we conclude that,
$$
\left(\frac{m^2(t)}{8\pi}+
\frac{\lm }{p+1} \left( \all+{t}\right)^{p+1}\mu(t)\right)^{'}
-\frac{1}{p+1}\left( \all+{t}\right)m^{'}(t)\leq 0,\mbox{ for a.a. } t\in [0, \thl].
$$

Therefore, we see that,
$$
-\frac{m^2(t)}{8\pi}-\frac{\lm}{p+1} \left( \all+{t}\right)^{p+1}\mu(t)+\frac{1}{p+1}\all m(t)
-\frac{1}{p+1}\int\limits_{t}^{\thl} m^{'}(s)s\, ds\leq 0,\; \forall\,t\in [0,\thl).
$$
Clearly, by using \rife{diffem}, we have that,
$$
-\int\limits_{t}^{\thl} m^{'}(s)s\, ds=tm(t)+\int\limits_{t}^{\thl} m(s)\,ds  =tm(t)+e(t),
$$
and we conclude that,
\beq\label{level1t}
-\frac{m^2(t)}{8\pi}-\frac{\lm}{p+1} \left( \all+{t}\right)^{p+1}\mu(t)+\frac{1}{p+1}\left(\all+{t}\right) m(t)
+\frac{1}{p+1}e(t)\leq 0,\; \forall\,t\in [0,\thl).
\eeq

Evaluating \rife{level1t} at $t=0$ we find that

$$
-\frac{m^2(0)}{8\pi}-\frac{\lm}{p+1} \all^{p+1}\mu(0)+\frac{1}{p+1}\all m(0)
+\frac{1}{p+1}e(0)=
$$
\beq\label{level2}
-\frac{\lm^2}{8\pi}-\frac{\lm}{p+1} \all^{p+1}+\frac{\lm}{p+1}\all
+\frac{2\lm^2}{p+1}\el
\leq 0,
\eeq

which is \rife{level0}. It is readily seen that the equality holds if and only if $\Gamma(t)$ is a disk for any $t$, whence if and only
if $\xil$ is radial and $\om=\mathbb{D}$. Here and in the sequel the radial symmetry is intended up to a translation. The inequality \rife{level01} is a straightforward consequence
of \rife{level0} and the fact that $\all\leq 1$. Concerning the characterization of the
equality sign in \rife{level01} we observe that if the equality holds, then necessarily $\all(1-\all^p)=0$ and in particular
the equality holds in \rife{level0}. Therefore, if the equality holds in \rife{level01}, then $\om=\mathbb{D}$ and either
$\all=0$ or $\all=1$. But if $\om=\mathbb{D}$ and $\all=1$ then $\el=E_0(\mathbb{D})=\frac{1}{16\pi}$, see \eqref{stima-energy}, and then $\el$ cannot be equal to $\frac{p+1}{16\pi}$ in this case. Therefore,
if the equality holds in \rife{level01}, then $\om=\mathbb{D}$ and $\all=0$. On the contrary, suppose that
$\om=\mathbb{D}$ and $\all=0$. Then, since $\om=\mathbb{D}$,
the equality holds in \rife{level0} and $\lm(\frac{p+1}{16\pi}-\el)=0$. But if $\lm=0$ then necessarily $\all=1$, and then
$\all$ cannot be zero in this case. Therefore, if $\om=\mathbb{D}$ and $\all=0$ then the equality holds in \rife{level01}.
This concludes the characterization of the equality sign in \rife{level01}.
\finedim

\bigskip

Next we prove the sharp positivity threshold.

\medskip

{\em Proof of Theorem \ref{thmvar}.}
We first prove that if $p>1$ and $(\lm,\pl)$ is any solution of $\prl$ with $\all=0$
then $\lm\geq \lm_{*}(\om,p)$ where the equality holds if and only if,
up to a translation, $\om=\mathbb{D}$. By defining,
$
R_{p+1}(w)=\dfrac{\ino |\nabla w|^2}{\left(\ino |w|^{p+1}\right)^{\frac{2}{p+1}}}$,
$w\in H^1_0(\om)\setminus\{0\}$,
then standard arguments in the calculus of variations and, since $\om$ is of class $C^3$,
elliptic regularity theory, show that $v$ is a classical $C^{2,r}_0(\ov{\om})$ solution of
\beq\label{eqvp}
-\Delta v=\mu v^p \mbox{ in }\om,\quad v=0 \mbox{ on }\pa\om,
\eeq
if and only if $\mu=\dfrac{R_{p+1}(v)}{\left(\ino |v|^{p+1}\right)^{\frac{p-1}{p+1}}}$
and $v$ is a critical point of $R_{p+1}$. In particular
$$
R_{p+1}(v)\geq \inf\limits_{w\in H^1_0(\om)}R_{p+1}(w)=\Lambda(\om,p+1),
$$
for any solution of \rife{eqvp}.
On the other side, if $(\lm,\pl)$ is any
solution of $\prl$ with $\all=0$, then for $\xil=\lm \pl$ solving \eqref{ul.1} we have $\el=\frac{1}{2\lm}\|\xil\|_{p+1}^{p+1}$ and
$$
\lm=\dfrac{R_{p+1}(\xil)}{\|\xil\|_{p+1}^{p-1}}\geq
\dfrac{\Lambda(\om,p+1)}{\|\xil\|_{p+1}^{p-1}}
=\dfrac{\Lambda(\om,p+1)}{\left(2\lm\el\right)^{\frac{p-1}{p+1}}}\,.
$$
As a consequence we readily deduce that,
$$
\lm^{2p}\geq \dfrac{\Lambda^{p+1}(\om,p+1)}{\left(2\el\right)^{{p-1}}},
$$
where the equality holds if and only if $\xil$ is a minimizer of $R_{p+1}$.
At this point, since by Theorem~\ref{energy} we have that $2\el\leq \frac{p+1}{8\pi}$,
where the equality holds if and only if, up to a translation, $\om=\mathbb{D}$, then
we also find that,

\beq\label{lastlms}
\lm^{2p}\geq \left(\frac{8\pi}{p+1}\right)^{p-1}{\Lambda^{p+1}(\om,p+1)},
\eeq
where the equality holds if and only if, up to a translation, $\om=\mathbb{D}$.
In other words, we see from \rife{lastlms} that if $(\lm,\pl)$ is a solution of $\prl$ with $\all=0$, then
$\lm\geq \lm_*(\om,p)$ and that the equality holds if and only if, up to a translation,
$\om=\mathbb{D}$. With this result at hand we can conclude the proof.\\
Clearly \rife{lastlms} holds also for $p=1$ since in this case the energy plays no role.
As a consequence we readily infer that if $\lm\leq \lm_*(\om,p)$, then either
$\all>0$ or $\all=0$ which is the case if and only if either $p=1$ and $\lm=\lm_*(\om,1)$
(since $\lm_*(\om,1)=\Lambda(\om,2)=\lm^{(1)}(\om)$) or $p>1$, $\lm=\lm_*(\mathbb{D},p)$ and $\om$ coincides up to a translation with $\mathbb{D}$. This proves the first part of the claim.\\
Next, by Theorem A,
we know that for variational solutions of $\prl$, $\all>0$ if and only if
$\lm<\lm^{**}(\om,p)$ and in particular that  if $\lm=\lm^{**}(\om,p)$ then $\all=0$.
As a consequence by the first part of the proof we have $\lm^{**}(\om,p)\geq \lm_*(\om,p)$ as well as
the characterization of the equality sign.
\finedim

\bigskip
\bigskip
\

We next consider the universal explicit $L^{\infty}$-bound. Here $\Gamma(p)$ is the Euler Gamma function.
\bpr\label{aprenrg} Let $p\in [1,+\ii)$ and $(\all, \pl)$ be a solution of {\rm $\prl$}. Then it holds,
\beq\label{green.17}
\|\pl\|_{\ii}\leq \frac{\widetilde{k}_p}{4\pi}\left(\all +{2\lm \el}\right)^{\frac{p}{p+1}}\,,
\eeq
where $\widetilde{k}_p=(\Gamma(p+2))^{\frac{1}{p+1}}$. In particular we have,
\beq\label{lmzero.0}
 \|\pl\|_{\ii}<\frac{ p+1}{4\pi}\left(1+\frac{\lm p}{8\pi}\right).
\eeq
\epr
\proof
Suppose first $\lm>0$. Let us define,
$$
k_p(\om)=\left(\ino G_{\sscp \om}^{p+1}(0,y)\right)^{\frac{1}{p+1}},
$$
then, by the Green representation formula and the Holder inequality we see that,
$$
\frac{\thl}{\lm}=\ino G_{\sscp \om}(0,y)\left(\all+{\xil(y)}\right)^p\leq
k_p(\om)
\left(\ino \left(\all+{\xil}\right)^{p+1}\right)^{\frac{p}{p+1}}=
$$
$$
k_p(\om)\left(\all +\ino \left(\all+{\xil}\right)^p\xil\right)^{\frac{p}{p+1}}=
k_p(\om)\left(\all +{2\lm\el}\right)^{\frac{p}{p+1}}.
$$
By a well know result in \cite{Ban} (or either by some estimates due to R. Talenti (\cite{Tal2})) one can see that,
$$
\ino G_{\sscp \om}^{p+1}(0,y)\leq \ino G_{\mathbb{D}}^{p+1}(0,y)=\frac{\Gamma(p+2)}{(4\pi)^{p+1}},
$$
and then we deduce that,
$$
\|\pl\|_{\ii}=\frac{\thl}{\lm}\leq \frac{\widetilde{k}_p}{4\pi}\left(\all +{2\lm\el}\right)^{\frac{p}{p+1}},
$$
which is \rife{green.17} for $\lm>0$. Letting now $\lm\to0^+$ and using the fact that $\pl$ depends continuously on $\lm$ (\cite{BJ2}) we conclude that \rife{green.17} holds for $\lm=0$ as well.\\

Concerning \rife{lmzero.0}, we deduce from \rife{green.17} and \rife{level01} that,
$$
\|\pl\|_{\ii}\leq \frac{ p}{4\pi}\frac{\widetilde{k}_p}{p}
\left(\all +\frac{\lm p}{8\pi}\frac{p+1}{p}\right)^{\frac{p}{p+1}}<
\frac{ p}{4\pi}\frac{1+p}{p}
\left(1 +\frac{\lm p}{4\pi}\frac{p+1}{2p}\right)^{\frac{p}{p+1}}\leq\frac{ p+1}{4\pi}\left(1+\frac{\lm p}{8\pi}\right),
$$
where we used $\widetilde{k}_p< 1+p$ and $(1+a)^\beta\leq1+\beta a$ for any $a\geq0, \beta\leq1$.

\finedim

\

We next turn to the estimates about the $\frac12$-$\frac1q$ thresholds. Recalling the definition of
$\ell(\om)$ in \eqref{l}, we start with the following preliminary result.
\ble\label{apriori2}
Let $p\in [1,+\ii)$ and $(\all, \pl)$ be a solution of {\rm $\prl$}. Then it holds,
\beq\label{nestIII.2}
\thl \leq \frac{\lm}{2\pi-\lm p}\,\ell(\om) \quad \mbox{for } \lm< \frac{2\pi}{p}.
\eeq
\ele
\proof

For $\lm=0$, \eqref{nestIII.2} is trivially satisfied. Consider now $\lm>0$. By the Green representation formula we have,
\beq\label{ngreen}
\frac{\thl}{\lm}=\frac{1}{2\pi} \ino G_0(y)\left(\all+\xil(y)\right)^p,
\eeq
where
$$
G_0(y)=2\pi G_{\om}(0,y).
$$

Now if $\all=1$ then $\xil\equiv 0$ and \eqref{nestIII.2} holds true. Thus we can assume w.l.o.g. that $\all\in [0,1)$. Therefore
$$
\left(\all+{\left.\xil\right|_{\pa \om}}\right)^p=\all^p<1
$$
and we define
$$
\om_{+}=\left\{y\in \om\,:\,\left(\all+{\xil(y)}\right)^p>1\right\}, \quad
\om_{-}=\left\{y\in \om\,:\,\left(\all+{\xil(y)}\right)^p\leq 1\right\},
$$
which both have nonempty interior (recall also \eqref{nm1}).\\
Since $G_0(y)>0$ in $\om$, by using \eqref{ngreen} and the inequality $ab\leq e^a-1+(b\log(b))\mathbbm{1}_{\{b\geq0\}}$, $a,\,b>0$, we find that,
$$
\frac{\thl}{\lm}<
\frac{1}{2\pi}\int\limits_{\om}\left(e^{G_0(y)}-1\right)+\frac{1}{2\pi}\int\limits_{\om_{+}}\left(\all+{\xil(y)}\right)^p
\log\left(\all+{\xil(y)}\right)^p.
$$

By a classical isoperimetric inequality due to Huber (\cite{Hub}) we have,

$$
\int\limits_{\om}e^{G_0(y)}\leq \frac{1}{2\pi}\left(\,\,\int\limits_{\pa\om}e^{\frac12 G_0(y)}\right)^2=\frac{1}{2\pi}|\pa\om|^2,
$$
and since $|\om|=1$ we conclude that,

$$
\int\limits_{\om}\left(e^{G_0(y)}-1\right)\leq \frac{1}{2\pi}|\pa\om|^2-1=\ell(\om).
$$
Therefore we find that,
$$
\frac{\thl}{\lm}< \frac{\ell(\om)}{2\pi} +\frac{1}{2\pi}\int\limits_{\om_{+}}\left(\all+{\xil(y)}\right)^p
\log\left(\all+{\xil(y)}\right)^p<
$$

$$
\frac{\ell(\om)}{2\pi}+\log\left(\all+{\thl}\right)^p
\frac{1}{2\pi}\int\limits_{\om_{+}}\left(\all+{\xil(y)}\right)^p<
\frac{\ell(\om)}{2\pi}+\log\left(\all+{\thl}\right)^p
\frac{1}{2\pi}\int\limits_{\om}\left(\all+{\xil(y)}\right)^p,
$$
that is
$$
\thl\leq  \frac{\lm  }{2\pi}\ell(\om)+
\frac{\lm  }{2\pi}\log \left(\all+{\thl}\right)^p<\frac{\lm  }{2\pi}\ell(\om)+\frac{\lm p }{2\pi}\thl,
$$
which, for $\lm p<{2\pi}$, immediately implies that \rife{nestIII.2} holds. \finedim

\bigskip
\bigskip

At this point we show a first bound from below for the boundary value $\all$.
\bpr\label{al1q}
Let $p\in [1,+\ii)$ and $(\all, \pl)$ be a solution of {\rm $\prl$}. Then it holds,
$$
\all> \frac1q \quad \mbox{for } \lm \leq \dfrac{4\pi}{e p\,\ell(\om)}\,.
$$
\epr
\proof
For $\lm=0$ we already know that $\all=1$ and the thesis holds true. We thus consider $\lm>0$. We argue by contradiction and assume that $\all\leq \frac1q$ for some $\lm\leq  \dfrac{4\pi}{ep\,\ell(\om)}$.\\
First of all, since $\ell(\om)\geq 1$, we have
\beq\label{a99}
\lm p\leq \dfrac{4\pi}{e\ell(\om)}\leq \frac{4\pi}{e}<2\pi, \quad \forall \lm\leq \dfrac{4\pi}{ep \,\ell(\om)}.
\eeq
Therefore we can use \rife{nestIII.2}, which yields,
\beq\label{4pie.1}
\thl\leq \frac1p\,\frac{2\ell(\om)}{e\ell(\om)-2}=:\frac1p\,a_0, \quad \forall \lm\leq \frac{4\pi}{ep\,\ell(\om)},
\eeq

and $a_0$ is always positive and well defined since $e\ell(\om)-2\geq e-2>0$.
We can assume w.l.o.g. that $a_0>1$, since otherwise we would find that
$$
\left(\all+{\thl}\right)^p\leq \left(\frac1q+{\thl}\right)^p \leq\left(1-\frac{1}{p}+\frac{a_0}{p}\right)^p\leq  1,
$$
which contradicts \rife{nm1}. Observe now that by using
$$
\left(\all+{\thl}\right)^p\leq \left(1+\frac{a_0-1}{p}\right)^p
$$
in \rife{ngreen} we have,
$$
\thl < \frac{\lm \kappa(\om)}{4\pi}\left(1+\frac{a_0-1}{p}\right)^p\leq
\frac{\lm }{4\pi}\left(1+\frac{a_0-1}{p}\right)^p,
$$
where
$$
\kappa(\om)=4\pi\ino G_\om(0,y)\leq4\pi \sup\limits_{x\in \om}\ino G(x,y)dy
$$
and we used a classical rearrangement result \cite{Tal2}, which speaks that,
$$
\kappa(\om)\leq4\pi \sup\limits_{x\in \om}\ino G(x,y)dy\leq
4\pi\sup\limits_{x\in \mathbb{D}}\int\limits_{\mathbb{D}}G_{\mathbb{D}}(x,y)=4\pi\int\limits_{\mathbb{D}}G_{\mathbb{D}}(0,y)= 1.
$$
As a consequence we conclude that,
\beq\label{nabove}
\mbox{if }\thl\leq \frac1p\,a_0,\mbox{ then }\thl <\frac1p\,\frac{1}{e}\left(1+\frac{a_0-1}{p}\right)^p,\quad\forall\,
\lm\leq \frac{4\pi}{ep\,\ell(\om)}.
\eeq

In view of \rife{nabove} we can iterate the argument and conclude in particular that, for any $n\geq 1$ such that $a_{n-1}>0$
it holds,
\beq\label{nabove.1}
\thl < \frac1p\,a_n=\frac1p\,h(a_{n-1}),\quad\forall\,\lm\leq \frac{4\pi}{ep\,\ell(\om)},
\eeq
where
$$
h(t)=\frac{1}{e}\left(1+\frac{t-1}{p}\right)^p, \;t\in (0,+\ii).
$$
However, it is trivial to check that if $t\in [1,3]$, then
$h(t)-t\leq \max\{h(1)-1,h(3)-3\}\leq \max\left\{\frac{1}{e}-1,\frac{1}{e}\left(1+\frac{2}{p}\right)^p-3\right\}<e-3$, for
any $p\geq 1$. Since
$$
1<a_0=\frac{2\ell(\om)}{e\ell(\om)-2}\leq\frac{2}{e-2}<3,
$$
then $a_1=h(a_0)\leq a_0+e-3$ and for any $n\geq 2$ such that $a_{n-1}> 1$, we have
$a_n\leq a_0+n(e-3)\leq \frac{2}{e-2}+n(e-3)$. Therefore $a_{n_1}\leq a_0+n_1(e-3)< 1$, for some $n_1\geq 2$,
and $\thl< \frac1p\,a_{n_1}< \frac1p$\,. As a consequence we conclude that,
$$
\left(\all+{\thl}\right)^p< \left(1-\frac{1}{p}+\frac{1}{p}\right)^p\leq  1,
$$
which contradicts once more \rife{nm1}. This is the desired contradiction which concludes the proof of Proposition \ref{al1q}.
\finedim

\

Finally, we derive the following universal explicit estimates about the $\frac12$ and the
positivity threshold.
\bpr\label{lemc12} Let $p\in [1,+\ii)$ and $(\all, \pl)$ be a solution of {\rm $\prl$}. Then it holds,
$$
\all> \frac12 \quad \mbox{for } \lm\leq \dfrac{4\pi}{ep}\,.
$$
Moreover,
\beq\label{alb2}
\mbox{if } \all=0 \quad \mbox{then } \lm > \frac{16\pi}{e (p+1)}\,.
\eeq
\epr
\brm\label{remlms} {\it Actually, we can prove that there exists an increasing function $g:[1,+\ii)\mapsto [\frac{8\pi}{e},+\ii)$
satisfying $g(t)\geq \frac{16\pi}{e (t+1)}$, $t\in [1,4]$,
$g(t)\geq \frac{16\pi}{e t}$, $t\in [4,16]$, $g(t)\geq \frac{16\pi}{e t}\frac{t+1}{t}$, $t\in [16,24]$,
$g(t)\geq \frac{24\pi}{e (t+1)}$, $t\in [24,48]$,
$g(t)\geq \frac{24\pi}{e t}\frac{3}{2}$, $t\in [48,+\ii]$, such that
if $(\all,\pl)$ is a solution of {\rm $\prl$} with $\all=0$ then $\lm>g(p)$. We skip the details of this fact which
can be derived by the same arguments used in the proof of Proposition \ref{lemc12}.}
\erm

\proof We postpone the proof of \eqref{alb2} and start to deduce the first assertion. By \eqref{alb2} we can actually assume that $\all>0$. Suppose then
by contradiction that $\all\leq \frac12$ for some $\lm\leq \dfrac{4\pi}{ep}$\,.
If
$$
1-\frac{16\pi \el}{p+1}\leq \all,
$$
then we deduce from \rife{level0} that $\frac{\lm }{8\pi}{(p+1)}\geq 1-\all^p$,
that is,
$$
\all^p\geq 1-\frac1e \frac{p+1}{2p}, \quad\mbox{whenever }\lm\leq \dfrac{4\pi}{ep}\,.
$$
Since $1-\frac1e \frac{p+1}{2p}\geq 1-\frac1e$, for any $p\geq 1$,
then we also have, $\frac12\geq \all\geq \all^p\geq 1-\frac1e$,
which is a contradiction. Therefore it holds,
$$
1-\frac{16\pi \el}{p+1}> \all,\; \mbox{ that is }\el< \frac{p+1}{16\pi}(1-\all).
$$
At this point we use \rife{green.17}, recalling $\xil=\lm \pl$, and deduce that,
$$
\thl\leq \frac{\lm  }{4\pi}\widetilde{k}_p\left(\all +{2\lm\el}\right)^{\frac{p}{p+1}}<
\frac{\lm }{4\pi}\widetilde{k}_p\left(\all +\frac{\lm p}{4\pi}\frac{p+1}{2p}(1-\all)\right)^{\frac{p}{p+1}}\leq
$$
\beq\label{c12}
\frac{\widetilde{k}_p}{e p}\left(\all +\frac{1}{e}\frac{p+1}{2p}(1-\all)\right)^{\frac{p}{p+1}},
\quad\mbox{whenever }\lm\leq \dfrac{4\pi}{ep}\,.
\eeq

The function $f(t,\al)=\frac{\widetilde{k}_t}{e}\left(\al +\frac{1}{e}\frac{t+1}{2t}(1-\al)\right)^{\frac{t}{t+1}}$,
$t\geq 1$, $\al\leq \frac12$, satisfies,
$$
f(t,\al)\leq f\left(t,\frac12\right)\leq \frac{t}{2}\leq t(1-\al), \forall\,t\geq 1, \forall \al\,\leq \frac12,
$$
and we readily infer from \rife{c12} that $\thl< (1-\all)$,
which is a contradiction to \rife{nm1}. This completes the proof of the first part of the claim.\\

We next turn to the estimate \eqref{alb2}. We first infer from
\rife{green.17} that, whenever $\all=0$, it holds,
$\thl\leq \frac{\lm }{4\pi}\widetilde{k}_p\left({2\lm\el}\right)^{\frac{p}{p+1}}$.
In particular, since by \rife{level01} we have $\el\leq \frac{p+1}{16\pi}$, then we conclude that,
\beq\label{green.18}
\thl\leq 2\widetilde{k}_p\frac{\lm }{8\pi}\left(\frac{\lm p}{8\pi}\frac{p+1}{p}\right)^{\frac{p}{p+1}},
\eeq
whenever $\all=0$. At this point we can prove \rife{alb2}. Assume by contradiction that for
some $\lm\leq \frac{16\pi}{e (p+1)}$ there exists a solution of $\prl$ with $\all=0$. Therefore,
after a straightforward evaluation, it follows from \rife{green.18} that
if $\lm\leq \frac{16\pi}{e p}\frac{p}{p+1}$ then,
$$
\thl\leq \frac{2\widetilde{k}_p}{p}\frac{2p}{e(p+1)}\left(\frac2e\right)^{\frac{p}{p+1}}.
$$
The function $f_1(t)=2{\widetilde{k}_t}\frac{2t}{e(t+1)}\left(\frac{2}{e}\right)^{\frac{t}{t+1}}$, $t\in [1,+\ii)$, satisfies,
$f_1(t)<t, \forall\,t\in [1,+\ii)$ and then we deduce that $\thl<1$ for $p\in [1,+\ii)$
which contradicts \rife{nm1} (with $\all=0$).
\finedim

\

We can now complete the proof of Theorem \ref{corplasma1}.\\

\emph{Proof of Theorem \ref{corplasma1}.}
The $L^{\infty}$ bound \eqref{est.intro} is proved in Proposition \ref{aprenrg},
while \eqref{pos2} follows by Propositions \ref{al1q} and \ref{lemc12}.
\finedim

\

\begin{center}
\textbf{Acknowledgments}
\end{center}

The authors are grateful to the anonymous referee for several insightful remarks.

\

\

\end{document}